\documentclass[11pt,reqno,tbtags]{amsart}
\usepackage{amssymb}
\usepackage{ifdraft}
\usepackage[square,numbers]{natbib}
\usepackage{color}
\numberwithin{equation}{section}

\allowdisplaybreaks

\newtheorem{theorem}{Theorem}[section]
\newtheorem{lemma}[theorem]{Lemma}

\newtheorem{corollary}[theorem]{Corollary}

\theoremstyle{definition}

\newtheorem{remark}[theorem]{Remark}

\theoremstyle{remark}

\newenvironment{romenumerate}[1][0pt]{% optional argument changes indentation
\addtolength{\leftmargini}{#1}\begin{enumerate}% gives (i), (ii) etc.
 }{\end{enumerate}}

\newcounter{oldenumi}
{\setcounter{oldenumi}{\value{enumi}}
\begin{romenumerate} \setcounter{enumi}{\value{oldenumi}}}
{\end{romenumerate}}

\newcounter{thmenumerate}

\newcounter{romxenumerate}   %less indented than standard.

\newcounter{xenumerate}   %no left indentation; thus wider lines

\newcommand{\refT}[1]{Theorem~\ref{#1}}
\newcommand{\refC}[1]{Corollary~\ref{#1}}
\newcommand{\refL}[1]{Lemma~\ref{#1}}
\newcommand{\refR}[1]{Remark~\ref{#1}}
\newcommand{\refS}[1]{Section~\ref{#1}}

\newcommand{\refand}[2]{\ref{#1} and~\ref{#2}}

\begingroup
  \divide\count255 by 60
  \multiply\count255 by -60
  \advance\count255 by \time
  \ifnum \count255 < 10 \xdef\klockan{\the\count1.0\the\count255}
  \else\xdef\klockan{\the\count1.\the\count255}\fi
\endgroup

\newcommand{\sumin}{\sum_{i=1}^n}
\newcommand{\sumij}{\sum_{i=1}^j}

\newcommand\set[1]{\ensuremath{\{#1\}}}

\newcommand\bigpar[1]{\bigl(#1\bigr)}
\newcommand\Bigpar[1]{\Bigl(#1\Bigr)}

\newcommand\lrpar[1]{\left(#1\right)}

\def\rompar(#1){\textup(#1\textup)}    % usage: \rompar(...)

\def\xexp(#1){e^{#1}}

\newcommand\floor[1]{\lfloor#1\rfloor}

\newcommand\ntoo{\ensuremath{{n\to\infty}}}

\newcommand\punkt[1]{\if.#1\else.\spacefactor1000\fi{#1}}
\newcommand\iid{i.i.d\punkt}    
\newcommand\ie{i.e\punkt}
\newcommand\eg{e.g\punkt}

\newcommand{\as}{a.s\punkt}

\newcommand{\tend}{\longrightarrow}
\newcommand\dto{\overset{\mathrm{d}}{\tend}}

\newcommand\eqd{\overset{\mathrm{d}}{=}}

\newcounter{CC}
\newcommand{\CC}{\stepcounter{CC}\CCx} %new constant C_i
\newcommand{\CCx}{C_{\arabic{CC}}}     %repeats the last C_i
\newcommand{\CCdef}[1]{\xdef#1{\CCx}}     %defines #1 as the last C_i
\newcounter{cc}
\newcommand{\cc}{\stepcounter{cc}\ccx} %new constant c_i
\newcommand{\ccx}{c_{\arabic{cc}}}     %repeats the last c_i
\newcommand{\ccdef}[1]{\xdef#1{\ccx}}     %defines #1 as the last c_i

\newcommand\E{\operatorname{\mathbb E{}}}
\renewcommand\P{\operatorname{\mathbb P{}}}
\newcommand\Var{\operatorname{Var}}

\newcommand\spann{\operatorname{span}}

\newcommand\gam{\gamma}
\newcommand\gG{\Gamma}

\newcommand\gs{\sigma}
\newcommand\gss{\sigma^2}

\newcommand\cE{\mathcal E}

\newcommand\cS{{\mathcal S}}
\newcommand\cT{{\mathcal T}}
\newcommand\cU{{\mathcal U}}

\newcommand\qw{^{-1}}

\newcommand\qq{^{1/2}}
\newcommand\qqw{^{-1/2}}
\newcommand\qqc{^{3/2}}
\newcommand\qqcw{^{-3/2}}

\renewcommand{\=}{:=}

\newcommand\intoi{\int_0^1}
\newcommand\intoo{\int_0^\infty}

\newcommand\dd{\,\textup{d}}

\newcommand\rv{random variable}

\newcommand\bfs{breadth first search}
\newcommand\rt{^{(t)}}
\newcommand\tS{\widetilde S}
\newcommand\GW{Galton--Watson}
\newcommand\GWtree{\GW{} tree}
\newcommand\cGWt{conditioned \GW{} tree}
\newcommand\GWp{\GW{} process}
\newcommand\txi{\widetilde\xi}
\newcommand\wx[1]{Z_{#1}}
\newcommand\wk{Z_k}
\newcommand\wl{Z_\elk}
\newcommand\wltn{\wl(\cT_n)}
\newcommand\wktn{\wk(\cT_n)}
\newcommand\ctn{\cT_n}
\newcommand\elk{k}
\newcommand\hxi{\hat\xi}
\newcommand\hhxi{\hat\xi'}
\newcommand\hT{\widehat\cT}
\newcommand\htx[1]{\hT^{(#1)}}
\newcommand\htk{\htx{k}}
\newcommand\pd{probability distribution}
\newcommand\tZ{\tilde Z}
\newcommand\REM[1]{{\raggedright\texttt{[#1]}\par\marginal{XXX}}}

\newenvironment{comment}{\setbox0=\vbox\bgroup}{\egroup} %deletes!

\newcommand\urladdrx[1]{{\urladdr{\def~{{\tiny$\sim$}}#1}}}

\begin{document}

\title[Sub-Gaussian tails for width and height of GW trees]
{Sub-Gaussian tail bounds for the width and height of conditioned Galton--Watson trees. } 

\date{17 November, 2010 
\ifdraft{\small (typeset \today{} \klockan)}
{}} 
%; revised ...

\author{Louigi Addario-Berry}
\address{Department of Mathematics and Statistics, McGill University, 805
  Sherbrooke Street West,  	Montr\'eal, Qu\'ebec, H3A 2K6, Canada}
\email{louigi@math.mcgill.ca}
\urladdrx{http://www.math.mcgill.ca/~louigi/}

\author{Luc Devroye}
\address{School of Computer Science, McGill University, 3480 University Street, 
Montr\'eal, Qu\'ebec, H3A 2A7, Canada}
\email{luc@cs.mcgill.ca}
\urladdrx{http://cg.scs.carleton.ca/~luc/}

\author{Svante Janson}
\address{Department of Mathematics, Uppsala University, PO Box 480,
 SE-751~06 Uppsala, Sweden}
\email{svante.janson@math.uu.se}
\urladdrx{http://www.math.uu.se/~svante/}

%\keywords{<keywords>}

\begin{abstract} 
We study the height and width of a Galton--Watson tree with offspring distribution 
$\xi$ satisfying $\E\xi=1$, $0 < \Var \xi < \infty$, conditioned on having 
exactly $n$ nodes. Under this conditioning, we derive 
sub-Gaussian tail bounds for both the width (largest number of nodes in any level) and height  
(greatest level containing a node); the bounds are optimal up to constant factors in the exponent. 
Under the same conditioning, we also derive essentially optimal upper tail bounds for 
the number of nodes at level $k$, for $1 \leq k \leq n$. 
\end{abstract}

\maketitle

\section{Introduction}\label{S:intro}

A \GWtree{} is the family tree of a \GWp, \ie, it is a random rooted tree,
constructed recursively from the root, where each node has a random number
of children and these random numbers are independent copies of some random
variable $\xi$ taking values in \set{0,1,\dots}. We let $\cT$ denote a
(random) \GWtree.
($\cT$ depends of course on $\xi$, or rather its distribution, but the
offspring distribution $\xi$ is fixed throughout the paper and is therefore
not shown explicitly in the notation.) We view the children of each node 
as arriving in some random order, so that $\cT$ is an ordered, or {\em plane} tree. 

At times in the paper it will be useful to think of $\cT$ as a subtree of 
the so-called {\em Ulam--Harris} tree $\cU$: this is the tree with root $\varnothing$ 
whose non-root nodes correspond to finite sequences of integers $v_1\ldots v_k$, 
with $v_1\ldots v_k$ having parent $v_1 \ldots v_{k-1}$ and children 
$\{v_1 \ldots v_{k} i~:~i \in \{1,2,\ldots\}\}$. For a node $v$ of $\cU$ we think of $vi$ as the 
$i$'th child of $v$. Any rooted plane tree $T$ in which all nodes have at most countably 
many children can be viewed as 
a subtree of $\cU$ by sending the root of $T$ to the root $\varnothing$ of $\cU$ 
and using the ordering of children in $T$ to recursively define an embedding of $T$ into $\cU$
(see \eg{} \cite{LeG}). 

We will study the \emph{\cGWt} $\ctn$, which is the random tree $\cT$
conditioned on having exactly $n$ nodes. 
In symbols,
$\ctn\=(\cT\mid|\cT|=n)$, where, for any tree $T$, $|T|$ denotes its number
of nodes. 
(We consider in the sequel only $n$ such that $\P(|\cT|=n)>0$.) 
For examples of standard types of random trees that can be represented as
\cGWt{s} for suitable $\xi$, see \eg{} \citet{Devroye}.
The \cGWt{s} are essentially the same as the %combinatiorially defined
random \emph{simply generated  trees} \cite{MM}, see \eg{}
\cite{Devroye} or \cite{Drmota}.

As is well-known, the distribution of the tree $\ctn$ is  not changed if
$\xi$ is replaced by another random variable $\xi'$ whose distribution is
replaced by \emph{tilting} (or \emph{conjugation})
\cite{Kennedy}:
$\P(\xi'=k)=c a^k\P(\xi=k)$, $k\ge0$, for some $a>0$ and normalizing
constant $c$. (Necessarily, $c=\bigpar{\E a^\xi}\qw$, and thus 
$\E a^\xi<\infty$.)
We may, except in some exceptional cases, by a suitable tilting assume that
$\E\xi=1$, so that the branching process is critical.
This turns out to be convenient, and 
we will in the sequel always make this assumption $\E\xi=1$. 
We further
assume that $\xi$ has finite variance $\gss\=\Var\xi<\infty$.
We exclude the trivial case $\xi=1$ a.s., \ie, we assume
$\gss>0$. (Equivalently, when $\E\xi=1$, $\P(\xi=0)>0$.)

For a rooted tree $T$ (deterministic or random), 
the \emph{depth} $h(v)$  of a node $v$ is its distance to the
root; the root thus has depth 0.
Let $\wl(T)$ be the
width at level $\elk$, \ie, the number of nodes at depth $\elk$,
$\elk=0,1,\dots$. We define, as usual, the \emph{width} of the 
tree by 
\begin{equation}
W=  W(T)\=\max_{\elk\ge0} \wl(T),
\end{equation}
and the \emph{height} by 
\begin{equation}
 H= H(T)\=\max\set{h(v):v\in T}=\max\set{\elk: \wl(T)>0}.
\end{equation}
%equivalently, $H(T)\=\max\set{h(v):v\in T}$.

It is well-known  that the width and height of a \cGWt{} $\ctn$ both are of
the order $\sqrt n$. 
More precisely, $n\qqw W(\ctn)$ and $n\qqw H(\ctn)$ both converge in
distribution, as \ntoo, see \eg{} \cite{AldousII}, \cite{ChMY},
\cite{DrmotaG:profile}  
and \cite{Drmota};
moreover, they converge jointly \cite{ChMY}, \cite{SJ167}, 
\begin{equation}\label{limres}
  \bigpar{n\qqw W(\ctn),n\qqw H(\ctn)}\dto(\gs W,\gs\qw H)
\end{equation}
for some limit variables $W$ and $H$,
that furthermore do not depend on the distribution of $\xi$.
($W$ is the maximum of a Brownian excursion, and $H\eqd 2W$;
see further \cite{SJ174}.)

Two of the main results of the paper are to prove essentially optimal uniform sub-Gaussian 
upper tail bounds for both $W(\ctn)/\sqrt{n}$ and $H(\ctn)/\sqrt{n}$ for every 
offspring distribution $\xi$ with finite variance. 
As an immediate consequence, the estimates $\E W(\ctn)=O(n\qq)$ and $\E H(\ctn)= O(n\qq)$ hold; 
even these much weaker statements are to our knowledge new at this level of
generality. 
(For estimates assuming an exponential moment of $\xi$, see \eg{} \cite{FGOR}.)

We let $C_1,C_2,\dots,c_1,c_2,\dots$ denote positive constants that may
depend on the distribution of $\xi$ (and in particular on $\gss$) but not on
$n$ or other parameters unless explicitly indicated.
(We use $C_i$ for ``large'' and $c_i$ for ``small'' constants.)
Proofs are given in \refS{Spf}.
\begin{theorem}
  \label{TW}
Suppose that $\E\xi=1$ and\/ $\Var\xi<\infty$. Then
\begin{equation*}
  \P\bigpar{W(\ctn)\ge x} \le \CC e^{-\cc x^2/n}
\CCdef\CCTW \ccdef\ccTW
\end{equation*}
for all $x\ge0$ and $n\ge1$.
\end{theorem}
\begin{theorem}\label{TH}
Suppose that $\E{\xi}=1$ and $0 < \Var \xi < \infty$. 
Then 
\begin{equation}\label{xsmhbig}
\P(H(\ctn)\geq h) \leq \CC\CCdef\CCsm e^{-\cc\ccdef\ccup h^2/n}
\end{equation}
for all $h \geq 0$ and $n \geq 1$. 
\end{theorem}
The condition $\Var \xi > 0$ excludes the case $\P(\xi=1)=1$, in which case
$\ctn$ is a path of length $n$.  
\begin{corollary}
    \label{CW}
Suppose that $\E\xi=1$ and\/ $0 < \Var\xi<\infty$. Then
$\E W(\ctn)=O(n\qq)$ and $\E H(\ctn)=O(n\qq)$. 
More generally, for every fixed $r<\infty$,
$\E(W(\ctn)^r)=O(n^{r/2})$ and $\E(H(\ctn)^r)=O(n^{r/2})$.
\end{corollary} 

While our methods do not prove the convergence \eqref{limres} of 
$W(\ctn)/\sqrt n$ and $H(\ctn)/\sqrt n$, we have thus as a corollary
obtained tightness of them, and we believe that our argument might be the
simplest proof of this tightness.

On the other hand, knowing the limit result \eqref{limres}, 
it follows from the fact that the bounds in \refC{CW} hold for every $r$ 
that all moments (also joint) converge in \eqref{limres}.
In particular, by the known formulas for the moments of $W$ and $H\eqd2W$
(see \eg{} \cite{BPY}),
as \ntoo,
\begin{align}\label{EW}
  \E\bigpar{W(\ctn)^r}/n^{r/2}
&\to \gs^r\E W^r = \gs^r2^{-r/2}r(r-1)\gG(r/2)\zeta(r), \\
  \E\bigpar{H(\ctn)^r}/n^{r/2}
&\to \gs^{-r}\E H^r = \gs^{-r}2^{r/2}r(r-1)\gG(r/2)\zeta(r). \label{EH}
\end{align}
For joint moments, see \cite{Donati-M} and \cite{SJ174}.
These results are well-known if $\xi$ is assumed to have an
exponential moment, see \eg{} \cite{flod} and \cite{DrmotaG:width}, 
but to our knowledge they have not, even in the case $r=1$, 
been proved before without extra conditions.

We emphasise that we obtain these bounds for higher moments of both 
$W(\ctn)$ and $H(\ctn)$, and even sub-Gaussian tail bounds for both variables,
without assuming more than a finite second 
moment of $\xi$. This is somewhat surprising, at least for the width, since
a $\xi$ with a large tail will produce a very wide Galton--Watson tree $\cT$
with comparatively large probability; the explanation is that if the tree
has one 
generation that is very large, say of size $m$,
then it will probably have many nodes (of order $m^2$) in later generations,
so the conditioning on exactly $n$ nodes makes this event very unlikely
if $m\gg\sqrt n$.
In other words, the bounds on the width hold, not because it is difficult
for the Galton--Watson tree to get many branches, but because it is
difficult to get rid of them in time.

\begin{remark} 
  We assume $\gss=\Var\xi<\infty$ throughout the paper. 
Since increasing $\gs$ makes 
the width larger and the height smaller (asymptotically at least), 
see \eg{} \eqref{EW}--\eqref{EH}, it is not reasonable
to expect that the results for the width generalize to the case
$\gss=\infty$.
However, for the same reason it seems likely that the results for the
height extend, but we have not investigated that and leave that as an open
problem. In particular, we ask the following questions (assuming $\E\xi=1$):
Is $\E H(\ctn)=O(n\qq)$ also if $\gss=\infty$?
Is $\E H(\ctn)=o(n\qq)$ if $\gss=\infty$?
\end{remark}

Next we consider the width $\wl(\ctn)$ at a given level $\elk$. 
%(In other words, the \emph{profile} of the tree $\cT_n$.)
Of course, $\wl(\ctn)\le W(\ctn)$, so the results above for $W(\ctn)$
immediately imply the same bounds for $\wl(\ctn)$, uniformly in $\elk$.
In particular,
\begin{equation}
  \label{ta}
\E \wl(\ctn)=O\bigpar{n\qq}.
\end{equation}
For $\elk\asymp n\qq$, this is the correct order of $\E\wltn$; in fact,
$n\qqw \wx{\floor{x\sqrt n}}(\ctn)$ converges in distribution for every
fixed $x\ge0$, and as a function of $x$,
see \cite{DrmotaG:profile,DrmotaG:width} 
(assuming a finite exponential moment)
and \cite{Kersting} (the general case, by probabilistic methods).

For small $\elk$, on the other hand, $\wltn$ is smaller and it was proven in
\cite[Theorem 1.13]{SJ167} that 
\begin{equation}\label{sofie}
  \E \wltn=O(\elk),
\end{equation}
uniformly for all $\elk\ge1$ and $n\ge1$.
This is the best possible estimate, since for any fixed $\elk$, 
\begin{equation}\label{sjw}
  \E \wltn\to1+\elk\gss,
\qquad\text{as \ntoo},
\end{equation}
see Meir and Moon \cite{MM} and Janson  \cite{SJ167,SJ188}. 
(It is shown in \cite{SJ188} that the sequence
  $\E \wltn$ is not always monotone in $n$, so \eqref{sofie} is not a
consequence of \eqref{sjw}.)

Furthermore, for large $\elk$, \eqref{sofie} is again not sharp. Indeed, if
$\elk\gg\sqrt n$, then typically $H(\ctn)<\elk$ and thus $\wltn=0$. In fact,
as $\elk\to\infty$, $\E\wltn$ decreases exponentially, as is shown by the
next theorem, which combines the three phases 
($\elk\ll\sqrt n$, $\elk\asymp n$, $\elk\gg\sqrt n$) in a unified statement.
(\citet{DrmotaG:width} gave
the  weaker bound $\CC n\qq e^{-\cc k/\sqrt n}$,
assuming an exponential moment on $\xi$.) 

\begin{theorem}
 \label{TWk}
Suppose that $\E\xi=1$ and $0<\Var\xi<\infty$.
For all $n,\elk\ge1$,
\begin{equation}\label{twk1}
  \E \wltn \le \CC \elk e^{-\cc \elk^2/n}
\end{equation}
and also
\begin{equation}\label{twk2}
  \E \wltn \le \CC  n\qq e^{-\cc \elk^2/n}
\CCdef\CCtwk \ccdef\cctwk
\end{equation}
(which is weaker for $\elk=o(\sqrt n)$ but equivalent for larger $n$).
\end{theorem}

Turning to higher moments of $\wltn$, we first note that for small $k$ there
is no result corresponding to \eqref{twk1} without assuming higher moments
of $\xi$. In fact, already for $k=1$, it is easy to see 
that for any $m\ge1$, 
$$\P(Z_1(\ctn)=m)\to m\P(\xi=m)$$ 
as \ntoo,
see \cite{Kennedy} and \refR{Rsizebiased}.
It follows by Fatou's lemma, that if $\E\xi^{r+1}=\infty$, for some $r>1$,
then $\E\wltn^r\to\infty$.
The same holds for $\E\wltn^r$ for every fixed $k\ge1$.

Conversely, it was proven in \cite[Theorem 1.13]{SJ167} that if 
$\E\xi^{r+1}<\infty$ for an integer $r\ge1$, then 
$\E\wltn^r =O(k^r)$ uniformly in $k\ge1$ and $n\ge1$.
(The restriction to integer $r$ is for technical reasons in the proof; we
conjecture that the result holds for any real $r\ge1$.)

On the other hand, the estimate \eqref{twk2} extends to higher moments
without assuming any moment condition on $\xi$ beyond our standing
$0 < \Var\xi<\infty$, \ie, $\E\xi^2<\infty$ and $\xi$ is not constant.

\begin{theorem}
 \label{TWkr}
Suppose that $\E\xi=1$ and $0<\Var\xi<\infty$.
For any $r<\infty$,
\begin{equation}\label{twkr}
  \E\bigpar{ \wltn /\sqrt n}^r\le \CC(r)  e^{-\cc \elk^2/n}
\end{equation}
for all  $\elk,n\ge1$.

Furthermore,
\begin{equation}\label{twke}
  \P\bigpar{ \wltn > x}\le \CC\CCdef\CCtwke  
e^{-\ccx\ccdef\cctwkek \elk^2/n-\cc\ccdef\cctwkex x^2/n}
\end{equation}
for all $x\ge 0$ and $n\ge1$.
\end{theorem}

\subsection{Remarks on the limit law.}
We say that $T$ is theta distributed if it has distribution function 
$$
\P(T\le x) 
= \sum_{j=-\infty}^\infty \left(1-2j^2x^2\right) e^{-j^2x^2} 
= \frac{4\pi^{5/2}}{x^3} \, \sum_{j=1}^\infty j^2 e^{-\pi^2 j^2 / x^2}, 
~ x > 0. 
$$
The appearance of $T$ as the limit law of the height of random conditional
Galton--Watson 
trees was 
noted in \cite{resz, deBKR, chung, k76, MM, flod}. 
Furthermore, the maximum of Brownian excursion of duration one is
distributed as $T/\sqrt{2}$ (see, e.g., \cite{BPY}).
In (\ref{limres}), $W \eqd T/\sqrt{2}$ and $H  \eqd T \sqrt{2}$.
It takes a moment to verify that for $x \ge 1$,
\begin{align}
\P(T \ge x)
%\ge 4x^2 e^{-x^2} 
&\ge  2  e^{-x^2},
\intertext{and for $x \le 1$,}
\P(T \le x)  
&\ge  40\, e^{-\pi^2 / x^2}. \label{tleft}  
\end{align}
The bound of \refT{TW}, combined with the limit result (\ref{limres})
then shows that
$$
\ccTW \le \frac{2}{\sigma^2}.
$$
Similarly, the bound of \refT{TH}, combined with the limit result (\ref{limres})
then shows that
$$
\ccup \le \frac{\sigma^2}{2}.
$$
It would be nice if $\ccTW$ and $\ccup$ could be be made more explicit.
In any case, the sub-Gaussian tail behaviour of the bounds in Theorems
\refand{TW}{TH}  
is optimal, modulo a constant factor (depending on $\xi$). 

We also have the trivial observation that
$$
W(\ctn) H(\ctn) \ge n-1.
$$
Thus, Theorems \refand{TW}{TH} yield the following left-tail upper bounds:
$$
\P(W(\ctn) \le x)
\le \P\left(H(\ctn) \ge \frac{n-1}{x} \right)
\le \CCsm \exp \left( - \frac{\ccup (n-2)}{x^2} \right)
$$
and
$$
\P (H(\ctn) \le x )
\le \P\left( W(\ctn)  \ge \frac{n-1}{x} \right)
\le \CCTW \exp \left( - \frac{\ccTW (n-2)}{x^2 } \right).
$$
In view of (\ref{limres}) and the remark \eqref{tleft}
about the theta distribution,
these bounds are optimal up to the constant factors $\ccTW$ and $\ccup$.

\section{Preliminaries}\label{Sprel}

The \emph{span} of $\xi$, denoted $\spann(\xi)$, is the largest integer $d$
such that $\xi/d$ \as{} is an integer. Note that $\P(|\cT|=n)>0$, so $\ctn$
exists, if and only if $n\equiv1$ modulo $\spann(\xi)$, except possibly for
some small $n$.

We let $\xi_i$ denote \iid{} copies of the random variable $\xi$, and let
$S_n$ be the partial sums of $\xi_1,\xi_2,\dots$,
\begin{equation}\label{sn}
  S_n\=\sumin\xi_i.
\end{equation}

By a classic formula, 
see \eg{}
Dwass \cite{Dwass}, Kolchin \cite[Lemma 2.1.3, p.~105]{Kolchin} or
\citet{Pitman:enum},
for $n\ge1$,
\begin{equation}\label{dwass}
  \P\bigpar{|\cT|=n}
=\frac 1n\P\bigpar{S_{n}=n-1},
\end{equation}
and, more generally, for 
$m,n\ge1$ and independent copies $\cT_1,\dots,\cT_m$ of $\cT$,
\begin{equation}\label{qk3}
  \P\Bigpar{\sum_{i=1}^m |\cT_i|=n}
=\frac mn\P\bigpar{S_{n}=n-m}.
\end{equation}
Together with the local central limit theorem, \eqref{dwass} implies
\cite[Lemma 2.1.4, p.~105]{Kolchin}, with $d\=\spann(\xi)$
(recall that we only consider $n$ such that $n\equiv 1\pmod d$), 
\begin{equation}\label{ptn}
  \P\bigpar{|\cT|=n}\sim \frac{d}{\sqrt{2\pi}\gs }\,n\qqcw.
\end{equation}

We will use a one-sided tail bound for $S_n$,
which we take from \citet{SJ167},
that only requires our (weak) conditions.
Note that, apart from the values of the constants, 
the bound in \eqref{l1} is exactly as the limit given by the local
central limit theorem when it applies; hence, at least for $m$ not too
large, it is of the best possible kind.

\begin{lemma}[{\cite[Lemma 2.1]{SJ167}}]
  \label{L1}
Suppose that $\xi_i$ are \iid, non-negative and integer-valued random
variables,  with $\E
\xi_i=1$ and $\Var\xi_i<\infty$, and let $S_n\=\sumin\xi_i$.
Then, for 
%some constants $\CC$ and $\cc$ (depending on the distribution of $\xi$) and 
all $n\ge1$ and $m\ge0$, 
\begin{equation}\label{l1}
  \P(S_n=n-m)\le \frac{\CCx}{\sqrt n} e^{-\ccx m^2/n}.
\CCdef\CCli \ccdef\ccli
\end{equation}
\end{lemma}

\begin{remark}
We can write the probability in \eqref{l1}
as $\P(\sumin(1-\xi_i)=m)$. The point is that
even without any assumptions on the tail of $\xi_i$ beyond finite variance,
the variables $1-\xi_i$ are bounded above, which is enough for strong tail
bounds for $m\ge0$. (There is no similar bound for $m<0$ under our weak
conditions.) Cf.\ the related tail bound
$\P(S_n\le n-m)\le{\CC} e^{-\cc m^2/n}$, which follows by \eqref{bern} below.
\end{remark}

We will use the following version of  Bernstein's inequality, 
which is valid for variables with a
one-sided bound, see \eg{} \cite[(2.9)--(2.13)]{Hoeffding} and 
\cite[Theorem  2.7]{McD}.

\begin{lemma}\label{Lbern}
  Let $X_1,X_2,\dots,X_n$ be independent random variables such that 
$X_i-\E X_i\le b$ for every $i$. 
Then, with $V\=\sum_{i=1}\Var(X_i) $,
  \begin{equation}\label{bern}
\P\lrpar{\sum_{i=1}^n (X_i-\E X_i) \ge t }	
\le \exp\lrpar{-\frac{t^2}{2 V+2bt/3}}.
  \end{equation}
\end{lemma}

\section{A size-biased \GWtree}\label{ShGW}

Let $\hxi$ be a \rv{} with the \emph{size-biased} distribution
\begin{equation}\label{hxi}
  \P(\hxi=m)=m\P(\xi=m).
\end{equation}
(Note that this is a \pd{} on \set{1,2,\dots} since $\E\xi=1$, and that
$\hxi\ge1$.) 

Let, for $k\ge1$, $\htk$ be the modified \GWtree{} defined as follows:
There are two types of nodes: \emph{normal} and \emph{mutant}. Normal nodes
have offspring 
(outdegree) according to independent copies of $\xi$, while mutant nodes
have offspring according to independent copies of $\hxi$. Moreover, 
all children of a normal node are normal, while
for each mutant node, one of its children is selected uniformly at random and 
called its \emph{heir}; the heir is mutant if it has depth less than $k$
but normal if the depth is at least $k$, and all other children are normal.
(Alternatively, we can call the mutants \emph{kings}, with a reproductive
behaviour different from the common people. At time $k$, a republic is
introduced, and everybody becomes equal.)

There are thus exactly $k$ mutant nodes, which together with the heir $v^*$
of the last mutant node form a path from the root to some node $v^*$ at
depth $k$. We call this path the \emph{spine} of $\htk$.

\begin{remark}\label{Rsizebiased}
  This construction with $k=\infty$ was introduced by \citet{LPP}, and is
  called the \emph{size-biased \GWtree}; in this case the spine is infinite
  so the tree is infinite. 
The underlying size-biased \GWp{} 
is  the same as the \emph{Q-process} studied in
\cite[Section I.14]{AN}.
For any fixed $k$, the first $k$ generations of $\ctn$ converge in
distribution to the first generations of $\htx\infty$.

Our $\htk$ is a truncated
  version of this, which grows like a normal \GWtree{} after generation $k$;
thus $\htk$ is \as{} finite.
\end{remark}

An equivalent construction is to start with the spine, and attach
independent copies of $\cT$ to it; the number of such trees attached to each
node in the spine except the last one (the top node)
has distribution $\hxi-1$, but the number
attached to the top node is $\xi$.

The probability that a given mutant node has $m$ children and that a given
one of them is selected as heir is, by \eqref{hxi},
\begin{equation*}
  \frac1m\P(\hxi=m)=\P(\xi=m), 
\qquad m\ge1.
\end{equation*}
It follows that for any rooted tree $T$,
and any path
$\gam$ in $T$ from the root to a node at depth $k$,
letting $d_1,d_2,\dots$ denote the outdegrees of the nodes in $T$, taken
in breadth-first order, say,
\begin{equation} \label{spineid}
  \P(\htk=T\text{ with $\gam$ as spine})
=\prod_v \P(\xi=d_v)
=\P(\cT=T).
\end{equation}
Since the possible spines in $T$ are in one-to-one correspondence with the
nodes at depth $k$, the number of them is $\wk(T)$, and thus
\begin{equation}\label{em}
  \P(\htk=T)=\wk(T)\P(\cT=T).
\end{equation}
In other words, $\htk$ has the distribution of $\cT$ biased by $\wk$, the
size of generation $k$. In particular, this yields, summing \eqref{em} over
all trees $T$ of size $|T|=n$,
\begin{equation*}
  \P\bigpar{|\htk|=n}=
\E\bigpar{\wk(\cT);\,|\cT|=n}
\end{equation*}
and thus
\begin{equation}\label{jb}
  \E\wktn=\frac{\E\bigpar{\wk(\cT);\,|\cT|=n}}{\P(|\cT|=n)}
=\frac{\P\bigpar{|\htk|=n}}{\P(|\cT|=n)}.
\end{equation}

\section{Proofs}\label{Spf}

\begin{proof}[Proof of \refT{TW}]
  Consider the \bfs{} of the Galton--Watson tree. As is well known,
this search keeps a queue of $Q_i$ nodes with $Q_0=1$ and the recursion
$Q_i=Q_{i-1}-1+\xi_i$, with $\xi_i$ \iid{} copies of $\xi$ as above;
hence $Q_j=1+\tS_j$, where $\tS_j\=\sum_{i=1}^j(\xi_i-1)=S_j-j$.
The \bfs{} stops, and the tree is completely explored, when
$Q_j$ becomes 0; in order for the tree to have size $n$ we thus have $Q_j>0$
for $0\le j<n$ and $Q_n=0$; equivalently, $\tS_j\ge0$ for $j<n$ and $\tS_n=-1$.

When the \bfs{} just has completed exploring the nodes at level $\elk-1$, the
queue consists of exactly the nodes at level $\elk$. Hence each $\wl$ is
some $Q_j$, and 
\begin{equation*}
W\=\max_{\elk\ge0} \wl\le \max_{j\ge0} Q_j.  
\end{equation*}
As a result, for the conditioned \GWtree{} $\ctn$, 
\begin{equation}
  \label{b2}
  \begin{split}
\P(W\ge x+1)
&\le \P\bigpar{\max_jQ_j\ge x+1}
\\&
=\P\bigpar{\max_j\tS_j\ge x\mid\tS_j\ge0,\, j<n, \text{ and } \tS_n=-1}.	
  \end{split}
\end{equation}
We get rid of the conditioning on $\tS_j\ge0$ ($j<n$) by the standard
rotation argument: for each (deterministic) sequence $x_1,\dots,x_n$ of
integers $\ge-1$ with sum $\sumin x_i=-1$, there is exactly one rotation
$x_i\rt\=x_{i+t}$ with $t\in\set{0,\dots,n-1}$ and indices taken modulo $n$,
such 
that the partial sums $S_j\rt\=\sumij x_i\rt\ge0$ for $1\le j< n$. Hence, we
can obtain $(\tS_j)_{j=1}^n$ with the conditional distribution given 
$\tS_j\ge0$, $j<n$, and  $\tS_n=-1$, as required in \eqref{b2}, by
conditioning $(\tS_j)_{j=1}^n$ on $\tS_n=-1$ and then taking the unique
correct rotation. The rotation may change $\max_j\tS_j$, but we have
$$
\max_{j\le n}\tS_j = \max_{j\le n}\tS_j - \min_{j\le n}\tS_j+1,
$$
and the latter quantity is changed by at most 1 by a rotation of 
$\txi_i\=\xi_i-1$, $i=1,\dots,n$. Hence, the rotation argument shows that
%\begin{equation}
%  \begin{split}
\begin{multline*}
\P\Bigpar{\max_{j\le n}\tS_j\ge x\mid\tS_j\ge0,\, j<n, \text{ and }
  \tS_n=-1}
\\
\le
\P\Bigpar{\max_{j\le n}\tS_j-\min_{j\le n}\tS_j
\ge x\mid\tS_n=-1}. 
\end{multline*}
%  \end{split}
%\end{equation}
By \eqref{b2} we thus have
\begin{multline*}
%\begin{equation}
%  \begin{split}
{\P(\max_j Q_j \ge 2x+2)} %\P(W\ge 2x+2)
\le
\P\bigpar{\max_{j\le n}\tS_j-\min_{j\le n}\tS_j \ge 2x+1\mid\tS_n=-1}
\\
\le
\P\bigpar{\max_{j\le n}\tS_j \ge x\mid\tS_n=-1}
+
\P\bigpar{\min_{j\le n}\tS_j \le -x-1\mid\tS_n=-1}.
%  \end{split}
%\end{equation}  
\end{multline*}
Furthermore, the reflection $\xi_i\leftrightarrow\xi_{n+1-i}$, which takes
$\tS_j\leftrightarrow\tS_n-\tS_{n-j}$, shows that the last probabilities are
the same, and we thus have
\begin{equation}\label{bem}
  \begin{split}
{\P(\max_j Q_j \ge 2x+2)} % \P(W\ge 2x+2)
\le
2\P\bigpar{\max_{j\le n}\tS_j \ge x\mid\tS_n=-1}.
  \end{split}
\end{equation}

Fix $x>0$ and let $\tau$ be the stopping time $\min\set{j\ge0:\tS_j\ge x}$.
Then \eqref{bem} can be written
\begin{equation}\label{b3}
  \begin{split}
{\P(\max_j Q_j \ge 2x+2)}% \P(W\ge 2x+2)
&\le
2\P\bigpar{\tau<n\mid\tS_n=-1}
\\&
=
\frac{2\P\bigpar{\tS_n=-1\mid\tau<n}\cdot\P(\tau<n)}
{\P(\tS_n=-1)}.
  \end{split}
\end{equation}
By definition, $\tS_\tau\ge x$.
Further, for any $t<n$ and $y\ge x$, by \refL{L1},
%\begin{equation*}  \begin{split}
\begin{multline*}
\P\bigpar{\tS_n=-1\mid\tau=t\text{ and }\tS_\tau=y}
=\P(\tS_n-\tS_t=-y-1)
\\
=\P\bigpar{\tS_{n-t}=-(y+1)}
\le \CCli n\qqw e^{-\ccli(y+1)^2/(n-t)} \le \CCli n\qqw e^{-\ccli x^2/n}.  
\end{multline*}
%  \end{split}\end{equation*}
Consequently,
\begin{equation*}
  \P\bigpar{\tS_n=-1\mid\tau<n}
\le \CCli n\qqw e^{-\ccli x^2/n},
\end{equation*}
and \eqref{b3} yields
\begin{equation}\label{maxq}
{\P(\max_j Q_j \ge 2x+2)} \le \frac{ \CC n\qqw e^{-\ccli x^2/n}}{\P(\tS_n=-1)}
\le \CC\CCdef\CCqf e^{-\cc\ccdef\ccqt x^2/n},
\end{equation}
since $\P(\tS_n=-1)\ge\cc n\qqw$ by the standard local central limit
theorem. 
Finally, since $\P(W \ge 2x+2) \le \P(\max_j Q_j \ge 2x+2)$, the proof is complete. 
%; recall that we assume $\spann\xi=1$.)
\end{proof}

\begin{proof}[Proof of \refT{TH}]
By choosing $\CCsm$ sufficiently large we may assume that $h \geq \sqrt{n}$. 
We may also assume that $h$ is an integer.  
Our proof of \eqref{xsmhbig} is based on the following observation: if $v$
is a node of $\ctn$ with ``large'' height then  
{\em either} there are many edges leaving the path from the root to $v$, {\em or} many of the ancestors of $v$ have exactly one child. 
In the first case, we will be forced to consider whether the majority of
edges leaving the root-to-$v$ path lead to nodes which are lexicographically
before, or after, $v$. To do so, we use {\em lexicographic} and 
{\em reverse-lexicographic} depth-first search (DFS) of $\ctn$.  

To define lexicographic DFS of $\ctn$, think of $\ctn$ as a plane tree (\ie{} as embedded in the Ulam--Harris tree $\cU$) and list the 
nodes of $\ctn$ in lexicographic order as $v_0,v_1,\ldots,v_{n-1}$. We then let 
$Q_0^d=1$ and $Q_{i}^d=Q_{i-1}^d-1+\xi_{v_{i-1}}$, where $\xi_{v_i}$ is the number of children of $v_i$ in $\ctn$. 
(This is sometimes called the Lukasiewicz path of $\ctn$; see, \eg{}, \cite{LeG}.) 
The {\em reverse-lexicographic} depth-first search of $\ctn$ is the sequence $Q_0^r,\ldots,Q_{|\ctn|}^r$ 
obtained by performing a lexicographic depth-first search on the mirror image of $\ctn$ (so if the root $\varnothing$ 
has children $1,\ldots,k$ in $\ctn$, then $k$ is the first rather than last child visited, and so on). We remark that the 
lexicographic and reverse-lexicographic depth-first search both
are identical in distribution to the breadth-first search of $\ctn$. 

Now let $p_1 = \P(\xi=1)$ and let $q_1=1-p_1$. If $v$ is a node of $\ctn$ with $h(v)=h$, then, writing $j$ (resp.~$k$) for the index of $v$ in lexicographic (resp.~reverse-lexicographic) order, 
either $\max(Q^d_j,Q_k^r)\geq (q_1/3)h$, or else at least $(p_1+q_1/3)h$ of the ancestors of $v$ have exactly one child. 
Let $\cS$ be the set of trees $T$ with $|T|=n$, such that $T$ contains a
node $v$ possessing $(p_1+q_1/3)h(v)$ 
ancestors with exactly one child and for which $h(v)=h$. Then let
$\cE\=\set{\ctn\in\cS}=\bigcup_{T \in \cS} \{\ctn=T\}$.  

Since $Q^d$ and $Q^r$ have the same distribution as $Q$, we then have 
\begin{align}
\P(H(\ctn) \geq h) & \leq \P(\max_j Q^d_j \geq (q_1/3)h) + \P(\max_k Q^r_k \geq (q_1/3) h) + \P(\cE) \nonumber\\
	& = 2\P(\max_i Q_i \geq (q_1/3)h) + \P(\cE) \nonumber\\
	& \leq \CC\CCdef\CCHf e^{-\cc\ccdef\ccHu h^2/n} + \P(\cE),\label{Hb}
\end{align}
the latter inequality holding due to (\ref{maxq}). 

Next, for each tree $T \in \cS$, fix a path $\gamma_T$ from the root of $T$ 
to a node $v$ with $h(v)=h$ and with at least $(p_1+q_1/3)h$ ancestors with
exactly one child  
(such a node exists by the definition of $\cS$). Then by 
(\ref{spineid}), 
\begin{equation}\label{tha}
  \begin{split}
%\P\left(\bigcup_{T \in \cS} \{\cT=T\}\right)& 
\P\left(\cT \in \cS\right)& 
= \sum_{T \in \cS} \P( \cT=T) \\
	& = \sum_{T \in \cS} \P(\hat{T}^{(h)}=T\text{ with $\gam_T$ as spine}) \\
& = \P\lrpar{\bigcup_{T \in \cS} \{\hat{T}^{(h)}=T\text{ with $\gam_T$ as spine}\}}\\
& \leq \P\left(\sum_{i=0}^{h-1} \mathbf{1}_{\hat{\xi}_i=1} \geq (p_1+q_1/3)h\right).	
  \end{split}
\end{equation}
The $\mathbf{1}_{\hat{\xi}_i=1}$ are Bernoulli$(p_1)$, so by \refL{Lbern}, 
\begin{align}
\P\left(\sum_{i=0}^{h-1} \mathbf{1}_{\hat{\xi}_i=1} \geq (p_1+q_1/3)h\right) 
& \leq \exp\left(-\frac{(q_1h/3)^2}{2p_1q_1h+2q_1^2h/9}\right) \nonumber\\
&= \exp\left(-\frac{h}{18p_1/q_1+2}\right). \label{thb}
\end{align}
It follows by \eqref{ptn} and \eqref{tha}--\eqref{thb} that 
\begin{align*}
\P(\cE) & = 
%\frac{\P\left(\bigcup_{T \in \cS} \{\cT=T\}\right)}{\P(|\cT|=n)}  \\
\frac{\P\left(\cT \in \cS\right)}{\P(|\cT|=n)}  
%	\\& \leq \frac{\sum_{T \in \cS} \P( \cT=T)}{\P(|\cT|=n)} 
	 \leq \CC n^{3/2} \exp\left(-\frac{h}{18p_1/q_1+2}\right)\\
	& \leq \CC\CCdef\CCat e^{-\cc\ccdef\ccau h^2/n}
\end{align*}
for all $h \geq \sqrt{n}$. 
Together with (\ref{Hb}) we have thus proved 
\[
\P(H(\ctn) \geq h) \leq \CCHf e^{-\ccHu h^2/n} + \CCat e^{-\ccau h^2/n},
\]
which establishes (\ref{xsmhbig}). 
\end{proof}

\begin{proof}[Proof of \refT{TWk}]
Note first that the case $k>n$ is trivial, since  $H(\ctn)\le n$ and
$\wktn=0$ for $k>n$. Further, if $k\le\sqrt n$, then the result follows from
\eqref{sofie}. Hence it suffices to consider $\sqrt n\le k\le n$.

Consider the random tree $\htk$ constructed in \refS{ShGW}. By the
alternative construction described there, we can regard the tree as the
$k$ mutant nodes (the spine except its top node) together with a random
number $M$ of attached independent copies of $\cT$. 
Hence, 
\begin{equation}\label{qk1}
|\htk|\eqd k+\sum_{i=1}^M |\cT_i|,
\end{equation}
where $\cT_i$ are independent copies of $\cT$, independent also of $M$.
The number $M$ is the
total number of normal children (including the top node) of the $k$ mutants, 
and thus 
\begin{equation}\label{qk2}
M\eqd\sum_{i=1}^k(\hxi_i-1)+1,   
\end{equation}
where $\hxi_i$ are \iid{} with the
distribution \eqref{hxi}.

Thus, for $m>0$ and $n> k$, using \eqref{qk1}, \eqref{qk3} and \refL{L1},
\begin{equation}\label{qk4}
  \begin{split}
\P\bigpar{|\htk|=n\mid M=m}
&=
\P\Bigpar{\sum_{i=1}^m |\cT_i|=n-k}
=
\frac{m}{n-k}\P\bigpar{S_{n-k}=n-k-m}
\\&
\le \CCli \frac{m}{(n-k)\qqc}e^{-\ccli m^2/(n-k)}.
  \end{split}
\raisetag\baselineskip
\end{equation}

The summands $\hxi_i-1$ in \eqref{qk2} have mean
$\E(\hxi-1)=\E\xi^2-1=\gss>0$.
We truncate them and define $\hhxi_i\=\min(\hxi_i,K)$, where $K$ is chosen so
large that $\E\hhxi_i>1+\gss/2$.
We apply Bernstein's inequality \eqref{bern} to $-\hhxi_i$, 
and obtain, since $\Var(\hhxi_i)<\infty$ and thus $V=O(n)$,
\begin{equation}\label{qk6}
  \begin{split}
\P\bigpar{M\le k\gss/4}
&\le \P\lrpar{\sum_{i=1}^k(\hxi_i-1)\le k\gss/4}	
\le \P\lrpar{\sum_{i=1}^k(\hhxi_i-1)\le k\gss/4}	
\\&
\le \P\lrpar{\sum_{i=1}^k(\hhxi_i-\E\hhxi_i)\le- k\gss/4}	
\le e^{-\cc k} \ccdef\ccb.
  \end{split}
\raisetag\baselineskip
\end{equation}

Note that $|\htk|\ge M+k$ by \eqref{qk1}, so if 
$M=m>k\gss/2$, we only have to consider $n\ge m+k>(1+\gss/2)k$, and for such
$n$, $n-k\ge \cc n$. Hence, for $m>k\gss/2$, \eqref{qk4} yields
\begin{equation}\label{qk5}
  \begin{split}
\P\bigpar{|\htk|=n\mid M=m}
\le \CC \frac{m}{n\qqc}e^{-\ccli m^2/n}
\le \CC \frac{1}{ n}e^{-\cc m^2/n}.
\CCdef\CCqkv \ccdef\ccqkv
  \end{split}
\end{equation}

If $\sqrt n\le k\le n$, \eqref{qk6} and \eqref{qk5} yield
\begin{equation}\label{magn}
  \begin{split}
\P\bigpar{|\htk|=n}
\le e^{-\ccb k}+\max_{m\ge k\gss/2}\CCqkv \frac{1}{n}e^{-\ccqkv m^2/n}
\le \CC \frac{1}{ n}e^{-\cc k^2/n}.
\ccdef\ccjesp
  \end{split}
\end{equation}

Since $\P(|\cT|=n)\ge\cc n\qqcw$ by \eqref{ptn}, 
\eqref{jb} and \eqref{magn} yield, if $\sqrt n\le k\le n$,
\begin{equation}
  \E\wktn \le \CC n\qq e^{-\ccjesp k^2/n} 
\le \CC k e^{-\ccjesp k^2/n} ,
%\ccdef\ccmagn
\end{equation}
which completes the proof. (We remarked above that it suffices to consider
such $k$.)
\end{proof}

\begin{proof}[Proof of \refT{TWkr}]
First, by \refT{TW},
\begin{equation*}
  \P(Z_k(\ctn)>x) \le \P\bigpar{W(\ctn)>x}
\le\CCTW e^{-\ccTW x^2/n}.
\end{equation*}
Further, 
since $Z_k(\ctn)>0$ implies $H(\ctn) \geq k$, 
\refT{TH} implies that 
\begin{equation*}
  \P(Z_k(\ctn)>x) \le \P\bigpar{H(\ctn)\ge k}
\le\CCsm e^{-\ccup k^2/n}.
\end{equation*}
Taking the geometric mean of these bounds we obtain
\eqref{twke}. 
Further, \eqref{twke}  implies, for any $r>0$,
with $\tZ\=Z_k(\ctn)/\sqrt n$,
\begin{equation*}
  \begin{split}
  \E \tZ^r &=r\intoo x^{r-1}\P(\tZ>x)\dd x
\le r\CCtwke\, e^{-\cctwkek k^2/n}\intoo x^{r-1}e^{-\cctwkex x^2}\dd x
\\&
=\CC(r) e^{-\cctwkek k^2/n}.
\qedhere
  \end{split}
\end{equation*}
\end{proof}

\section*{Acknowledgements} 
Part of this work was completed while 
all three authors were visiting McGill's Bellairs 
Research Institute, Barbados,
and another part while S.J.~was visiting the 
Centre de Recherches Math\'ematiques in Montr\'eal.

\newcommand\AAP{\emph{Adv. Appl. Probab.} }
\newcommand\JAP{\emph{Journal of Applied Probability} }
\newcommand\JAMS{\emph{J. \AMS} }
\newcommand\MAMS{\emph{Memoirs \AMS} }
\newcommand\PAMS{\emph{Proc. \AMS} }
\newcommand\TAMS{\emph{Trans. \AMS} }
\newcommand\AnnMS{\emph{Ann. Math. Statist.} }
\newcommand\AnnPr{\emph{Ann. Probab.} }
\newcommand\CPC{\emph{Combin. Probab. Comput.} }
\newcommand\JMAA{\emph{J. Math. Anal. Appl.} }
\newcommand\RSA{\emph{Random Structures and Algorithms} }
\newcommand\ZW{\emph{Z. Wahrsch. Verw. Gebiete} }
\newcommand\DMTCS{\jour{Discr. Math. Theor. Comput. Sci.} }

\newcommand\AMS{Amer. Math. Soc.}
\newcommand\Springer{Springer-Verlag}
\newcommand\Wiley{Wiley}

\newcommand\vol{\textbf}
\newcommand\jour{\emph}
\newcommand\book{\emph}
\newcommand\inbook{\emph}
\def\no#1#2,{\unskip#2, no. #1,} %(typeset after year) 
\newcommand\toappear{\unskip, to appear}

\newcommand\webcite[1]{\hfil  %???
   \penalty0 %???
\texttt{\def~{{\tiny$\sim$}}#1}\hfill\hfill}
\newcommand\webcitesvante{\webcite{http://www.math.uu.se/~svante/papers/}}
\newcommand\arxiv[1]{\webcite{arXiv:#1.}}

\end{document}